\documentclass{amsart}
\usepackage{amsmath,amssymb,amsthm}
\newtheorem{thr}{Theorem}[section]

\newtheorem{conj}[thr]{Conjecture}

\theoremstyle{definition}

\theoremstyle{remark}

\numberwithin{equation}{section}

\begin{document}


\title{The ring $\mathrm{M}_{8k+4}(\mathbb{Z}_2)$ is nil-clean of index four}

\author{Yaroslav Shitov}

\email{yaroslav-shitov@yandex.ru}


\begin{abstract}
We show that the direct sum of an \textit{odd} number of matrices 
$$C=\left(\begin{array}{cccc}
0&0&0&1\\
1&0&0&0\\
0&1&0&0\\
0&0&1&1
\end{array}\right)$$
cannot be a sum $P+Q$ of matrices over $\mathbb{F}_2$ satisfying $P^2=P$ and $Q^3=O$.
\end{abstract}

\maketitle

A ring $R$ is called \textit{nil-clean of index $k$} if every element $r\in R$ can be written as $r=p+q$ with $p^2=p$ and $q^k=0$, and $k$ is the smallest integer with this property. J.~\v{S}ter~\cite{Ster} proved that $\mathrm{M}_n(\mathbb{Z}_2)$ is a nil-clean ring of index at most $4$. For $n=4$ and $n=8$, he showed that the index is actually equal to four using the matrices $C$ and $C\oplus C$, respectively. Is the index equal to $3$ for all sufficiently large $n$? This question originates from~\cite{Ster} and also appears in~\cite{Breaz}; we solve it in the negative.


\begin{proof}
Let $P,Q$ be matrices as in the abstract; from $P^2=P$, $Q^3=O$, and since $t^4+t^3+1$ is an annihilating polynomial of $C$, it follows after simplifications that
\begin{equation}\label{eq1}
PQPQ+QPQP+PQ^2P+Q^2PQ+QPQ^2+PQP=I.
\end{equation}
We write $\alpha:=\operatorname{corank}P$, and an appropriate basis change leads us to
$$P=\left(\begin{array}{c|c}
I&O\\\hline
O&O
\end{array}\right),\,\,\,\,
Q=\left(\begin{array}{c|c}
Q_1&Q_2\\\hline
Q_3&Q_4
\end{array}\right)$$
with bottom-right blocks of both displayed matrices being $\alpha\times\alpha$. Passing to these blocks again in~\eqref{eq1}, we get $(Q_3Q_2)Q_4+Q_4(Q_3Q_2)=I_{\alpha}$ and take the trace to see that $\alpha$ is even. So $\operatorname{trace}(P+Q)=\operatorname{trace}(P)$ is zero, which cannot happen because $P+Q$ is the direct sum of an \textit{odd} number of copies of $C$.
\end{proof}


Can $C$ be a block of the rational form of the sum of any projector and an index-three nilpotent? The results of~\cite{Ster} and this paper suggest a negative answer.


\smallskip


\begin{conj}
The ring $\mathrm{M}_n(\mathbb{Z}_2)$ is nil-clean of index $4$ for all $n\geqslant 8$.
\end{conj}




\end{document}